\input amstex
\input xy
\xyoption{all}
\loadeusm

\documentstyle{amsppt}
\magnification=\magstep{1}  %Setting font to "12 pt"%
\pagewidth{6.5truein}
\pageheight{9.0truein}

%%\long\def\ignore#1\endignore{INSERT DIAGRAM}
\long\def\ignore#1\endignore{#1}

\def\edge{\ar@{-}}

\def\la{\Lambda}
\def\lamod{\Lambda{}\operatorname{-mod}}
\def\Lamod{\Lambda{}\operatorname{-Mod}}
\def\pdim{\operatorname{p\,dim}}
\def\gldim{\operatorname{gl\,dim}}

\def\findim{\operatorname{fin\,dim}}
\def\Findim{\operatorname{Fin\,dim}}
\def\ker{\operatorname{ker}}

\def\soc{\operatorname{soc}}
\def\be{{\bold e}}
\def\NN{{\Bbb N}}

\def\Q{{\Cal Q}}

\def\Bass{{\bf 1}}
\def\But{{\bf 2}}
\def\Fos{{\bf 3}}
\def\pre{{\bf 4}}
\def\dom{{\bf 5}}
\def\dep{{\bf 6}}
\def\JenLen{{\bf 7}}
\def\Scho{{\bf 8}}
\def\Sma{{\bf 9}}

\topmatter

\title Stacks of algebras and their homology \endtitle

\rightheadtext{Stacks of algebras}

\author Nancy Heinschel and Birge Huisgen-Zimmermann \endauthor

\address Department of Mathematics, University of California, Santa
Barbara, CA 93106, USA\endaddress

\email heins\@math.ucsb.edu\endemail

\address Department of Mathematics, University of California, Santa
Barbara, CA 93106, USA\endaddress

\email birge\@math.ucsb.edu\endemail

\thanks The first author was supported in part by a fellowship
stipend from the National Physical Science Consortium and the
National Security Agency. The second author was partially
supported by a grant from the National Science Foundation.
\endthanks

\dedicatory Dedicated to Raymundo Bautista and Roberto Martinez-V\'illa
on the occasion of their sixtieth birthdays \enddedicatory

\abstract  For any increasing function $f: \NN
\rightarrow \NN_{\ge 2}$ which takes only finitely many distinct
values, a connected finite dimensional algebra $\la$ is constructed, with
the property that $\findim_n \la = f(n)$ for all $n$;  here $\findim_n
\la$ is the $n$-generated finitistic dimension of $\la$.  The
stacking technique developed for this construction of homological
examples permits strong control over the higher syzygies of
$\la$-modules in terms of the algebras serving as layers.
\endabstract    

\endtopmatter

\document

\head  1. Introduction and background \endhead

The purpose of this paper is twofold.  One of our objectives is to
introduce a technique of `stacking' finite dimensional algebras on top of
one another so that, on one hand, the homology of the resulting algebra
can be controlled in terms of the layers, while, on the other hand, this
homology differs qualitatively from that of the building blocks.   Our
second, principal, goal is to apply such stacks towards realizing new
homological phenomena.

Given a finite dimensional algebra $\la$ over a field $K$ and $n \in \NN$,
we denote by
$\findim_n \la$ the supremum of the finite projective dimensions attained
on left
$\la$-modules with `top multiplicities $\le n$'; in other words,
if $J$ denotes the Jacobson  radical of $\la$, we are focusing on those
left $\la$-modules $M$ of finite projective dimension for which the
multiplicities of the simple summands of $M/JM$ are bounded above by
$n$.  Since, over a basic algebra, this condition just means that
$M$ can be generated by
$\le n$ elements, we refer to $\findim_n \la$ as the (left) {\it
$n$-generated finitistic dimension\/} of
$\la$.  While it is still open whether the little finitistic dimension,
$$\findim \la = \sup_{n \in \NN}\, \findim_n \la,$$ 
is always finite  -- 
the question goes back to Bass's 1960 paper \cite{\Bass}  --   it is
well known that $\findim_n
\la < \infty$ for all $n$ (\cite{\JenLen,  Proposition 10.33} and
\cite{\Scho}).  This puts a spotlight on jumps $\findim_n \la < \findim_m
\la$ for $n < m$.  Indeed, producing a counterexample to finiteness of
the little finitistic dimension would amount to constructing a finite
dimensional algebra
$\la$ together with an infinite sequence $n_1 < n_2 < n_3 < \cdots$ of
positive integers such that
$$\findim_{n_k} \la < \findim_{n_{k+1}} \la$$
 for all $k$.  While so far
there has been hardly any insight into the mechanism of such jumps  -- 
the first illustrations were based on monomial algebras, where only a
single jump is possible \cite{\dep} 
--  we now use our stacking technique to systematically create examples of
`compounded jumps of arbitrary size'.  Namely, for any 
increasing function $f: \NN \rightarrow
\NN_{\ge 2}$ which takes only finitely many distinct values, we
exhibit a connected finite dimensional algebra $\la$ with vanishing
radical cube such that 
$\findim_n \la = f(n)$ for all $n \in \NN$ (Section 4).  The recursion
underlying these successively `stacked' examples pinpoints a typical
combinatorial pattern giving rise to skips in the 
finitistic dimensions attained on modules of bounded length. 

In this connection, we note that arbitrary finite gaps
$\Findim \la - \findim \la$  between the little and big finitistic
dimensions are already known to be realizable in the finite dimensional
setting; here the big finitistic dimension, $\Findim \la$, is the
supremum of all finite projective dimensions attained on {\it
arbitrary\/} left
$\la$-modules.  Gaps of $1$ were first obtained for monomial
algebras in \cite{\dom}, where, however, they cannot exceed $1$;  see
\cite{\pre}.  Arbitrary jumps were produced by 
Rickard (unpublished), who showed that the big and little finitistic
dimensions are additive on tensor products,
and by Smal\o\, \cite{\Sma}, who constructed simpler examples via
iterated one-point extensions of an algebra
$\la_0$ with $\Findim \la_0 - \findim \la_0 = 1$.  

The stacking technique
presented in Section 3 below is governed by less restrictive rules than
one-point extensions:   Any one-point extension of an algebra
$\la_0$ over a field $K$ results from stacking the one-dimensional
$K$-algebra on top of $\la_0$, but stacking an algebra
$\la_1$ on top of an algebra $\la_0$ will not even lead to an {\it
iterated\/} one-point extension of $\la_0$ in general.  On the other hand,
this process always results in a triangular matrix algebra
$$\la = \pmatrix \la_1 & 0 \\ M & \la_0 \endpmatrix$$ 
for a suitable $\la_0$-$\la_1$-bimodule $M$.  Hence initial upper and
lower bounds on $\Findim \la$ in terms of $\la_0$ and $\la_1$ are
available: Indeed, Fossum, Griffith and Reiten show in
\cite{\Fos, Corollary 4.21} that  
$$\Findim \la_0 \le \Findim \la \le \Findim \la_0 + \Findim \la_1 + 1;$$
the little finitistic dimensions are subject to analogous
inequalities, as can easily be seen by the same method.  Our
stacks are considerably more specialized than the triangular matrix
construction, so as to allow for tighter control of syzygies over the new
algebra. 

Since the algebras we target in our examples are stacks of monomial
algebras, i\.e., of algebras of the form
$K\Q/I$, based on a field
$K$, a quiver $\Q$ and an admissible ideal $I$ that can be generated by
paths in 
$K\Q$, we include a brief review of the homology of monomial algebras at
the end of this section for easy reference.  These algebras are
homologically well-understood.  In particular, all second syzygies of
their modules are direct sums of cyclic left ideals recruited from a
finite collection (see
\cite{\dom} and Theorem 1 below);  moreover, the projective dimensions of
these cyclic left ideals can easily be computed (for an algorithm, see
\cite{\pre}).  The idea of our applications is to `stack' the complexity
observed in first syzygies over monomial algebras, while
benefiting from the simplicity of second syzygies.    

In Section 2, we will briefly discuss the graphs which we use to
communicate certain types of modules in an intuitive
format.  
\bigskip  

Throughout,  $\Lambda$ will be a split basic finite dimensional algebra
over an arbitrary field $K$.  The category of all left
$\la$-modules will be denoted by $\Lamod$, and $\lamod$ will be the full
subcategory having as objects the finitely generated modules.  For
simplicity, we will identify
$\la$ with a path algbra modulo relations,
$K\Q/I$, where
$\Q$ is a quiver and $I \subseteq K\Q$ an admissible ideal in the path
algebra.  If $p$ and $q$ are paths in $\Q$, then $pq$ will stand for `$q$
followed by $p$'.   Moreover, we will identify the set $\Q_0$ of vertices
of $\Q$ with a full set of primitive idempotents of $\la$, loosely
referred to as `the' primitive idempotents of $\la$.  

Given any (left) $\la$-module $M$, we call $x \in M$ a {\it top
element\/} of $M$ in case $x \in M \setminus JM$ and $x = ex$ for some
primitive idempotent $e$;  in this situation we also say that $x$ is a
{\it top element of type\/} $e$ of $M$.   For $e \in \Q_0$, the simple
module $\la e/ Je$ will be denoted by
$S(e)$.

All of the finitistic dimensions, $\findim_n \la$,
$\findim \la$, and $\Findim \la$, depend on the side.  However, there is
no need to weigh down our notation with left-right qualifiers, since we
will consistently deal with {\it left\/} modules.      
\medskip

For the remainder of the introductory section, we assume $\la$ to be a
monomial algebra.  This means that the set of
 paths in $K\Q\setminus I$ gives rise to a $K$-linearly independent
family of residue classes in $\la$, which will be called the {\it
nontrivial paths\/} in
$\la$; thus the nontrivial paths in $\la$ form a basis for $\la$
over $K$.  Clearly, it makes sense to speak of the {\it length\/} of a
nontrivial path in $\la$.   The paths in the following set  will be
called the {\it critical paths\/}:
$$\multline 
\eusm{P} = \{ p \in \Lambda \mid p \text{ is a nontrivial path of
positive length,} \\ 
\text{starting in a non-source of\ } \Q, \text{\ with p
dim}_{\Lambda}\Lambda p < \infty \}.
\endmultline $$   We use this set to define an invariant $\bold s$ as
follows:
$$ \bold s = \cases -1 &\text{ if } \eusm{P} = \emptyset \\
 \text{max} \{ \text{p dim}_{\Lambda}\Lambda p \,|\, p \in
\eusm{P} \}  &\text{ otherwise.} \endcases$$  
The set $\eusm{P}$, as well
as the number $\bold s$, can be readily obtained from the graphs of the
indecomposable left
$\Lambda$-modules. 

On one hand, the smallest in our gamut of finitistic dimensions,
$\findim_1 \la$, is trivially bounded below by $\bold{s} + 1$. On
the other hand, the following theorem (see \cite{\dom}) shows that all
finitistic dimensions are bounded above by
$\bold{s} + 2$.  A generalization of the first part, together with a slick
argument, can be found in \cite{\But}.
 
\proclaim{Theorem 1} Let $M$ be a submodule of a projective left
$\Lambda$-module, and $E(M)$ the set of those primitive idempotents $e$ of
$\la$ which do not annihilate $M/JM$. 

\noindent \rom{(1)} The syzygy $\Omega _{\Lambda}^1(M)$ of $M$ is
isomorphic to a direct sum of principal left ideals of $\Lambda$, each
generated by a nontrivial path of positive length starting in a vertex in
$E(M)$.  

In particular, all second syzygies of $\la$-modules are direct sums of
cyclic left ideals $\la p$, $p$ a path, and all finitistic dimensions of
$\la$  fall into the interval $[\bold{s} + 1,\bold{s} + 2]$.
\smallskip  
 
\noindent \rom{(2)} 
Given any nontrivial path $q$ of positive length in $\Lambda$, the
following statements are equivalent: 
\roster
\item"(i)" $\Lambda q$ is isomorphic to a direct summand of $\Omega
_{\Lambda}^1(M)$.
\item"(ii)" There exists a path of positive length in $\la$, say $\alpha
p$ where $\alpha$ is an arrow, together with a top element $x$ of $M$
such that
$\Lambda q \cong \Lambda \alpha p$ and $px
\notin J^{\text{length}(p) + 1}M$, while $\alpha p x = 0.$
\endroster
\endproclaim 

From the first part of Theorem 1 we glean that, over a monomial algebra
$\la$, there is at most one positive integer $n$ with the property that 
$\findim_n \la <  \findim_{n+1} \la$, and these two dimensions
differ by at most $1$.

\head 2.  Graphs of modules \endhead

The modules arising in our examples can be represented by {\it layered}
undirected graphs of a format which is intuitively suggestive.  We use
the conventions from \cite{\pre, Section 5} and \cite{\dep, Section 2};
but for the reader's convenience, we briefly review the very simple
special cases needed here.

Let $\la = K\Q/I$ be a path algebra modulo relations  --  it will
reappear as $\la_1$ in Section 4  -- with quiver $\Q$ as follows.
$$\xymatrixrowsep{3.0pc} 
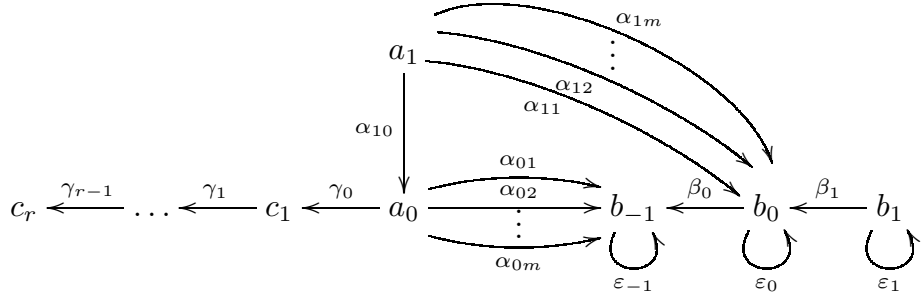
\xymatrix{  & & & a_1 \ar[d]_{\alpha _{10}}
\ar@(r,ul)[drrr]<-0.7ex>_(0.4){\alpha _{11}}
\ar@(r,ul)[drrr]<1.6ex>_(0.45){\alpha _{12}}^{\vdots}
\ar@(ur,u)[drrr]<1.7ex>^{\alpha _{1m}}\\
  c_r & \dots \ar[l]_{\gamma _{r-1}} & c_1 \ar[l]_{\gamma _1} & a_0
\ar[l]_{\gamma _0} \ar@/^/[rr]<0.85ex>^{\alpha _{01}} \ar[rr]^{\alpha
_{02}} \ar@/_/[rr]<-1.5ex>^{\vdots}_{\alpha _{0m}} & & b_{-1}
\ar@(dl,dr)_{\varepsilon _{-1}} & b_0 \ar[l]_{\beta _0}
\ar@(dl,dr)_{\varepsilon _0} & b_1 \ar[l]_{\beta _1}
\ar@(dl,dr)_{\varepsilon _1}
 }$$ 

\noindent A generating set for the ideal $I$ of relations can be
communicated by way of graphs of the indecomposable projective left
$\la$-modules, as presented in Part I of the proof of Theorem 10. 
We give a few samples to explain how to interpret these graphs.  While,
in general, our graphs are to be read relative to a given sequence of top
elements
$x_i$ of the considered module $M$ (generating $M$ modulo $JM$ and
linearly independent modulo
$JM$), when $M =
\la e$ is indecomposable projective, we tacitly assume the choice of top
element to be
$x = e$.

Thus, presenting $M = \la a_0$ by way of the layered graph
$$\xymatrixcolsep{0.75pc} \xymatrixrowsep{1.0pc}
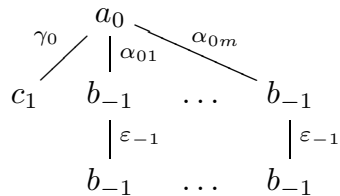
\xymatrix{ & a_0 \edge[dl]_{\gamma _0} \edge[d]^{\alpha _{0 1}}
\edge[drr]^{\alpha _{0 m}} \\ c_1 & b_{-1} \edge[d]^{\varepsilon
_{-1}} & \dots & b_{-1} \edge[d]^{\varepsilon _{-1}} \\
 & b_{-1}& \dots &
b_{-1}  }$$ 
holds the following information: The left ideal $Ia_0 \subseteq K\Q$ is
generated by
$\gamma_1\gamma_0$ and $\varepsilon_{-1}^2 \alpha_{0i}$ for $1 \le i \le
m$.  Equivalently, $M/JM \cong S(a_0)$,
$JM/J^2M \cong S(c_1) \oplus S(b_{-1})^m$  and $J^2M \cong S(b_{-1})^m$,
while $J^3 M = 0$. 

That $M = \la a_1$ has graph
$$\xymatrixcolsep{1.5pc}
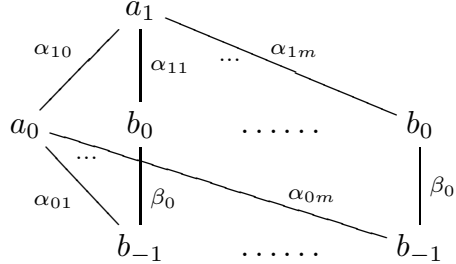
\xymatrix{ & a_1 \edge[dl]_{\alpha_{1 0}} \edge[d]^{\alpha _{1 1}}
\edge[drr]^{\alpha _{1 m}}_(0.35){\dots}  \\ a_0 \edge[dr]_{\alpha_{0
1}}^(0.35){\dots} \edge[drrr]^(0.7){\alpha _{0 m}} & b_0
\edge[d]^(0.6){\beta _0} & \dots\dots & b_0 \edge[d]^{\beta _0} \\ &
b_{-1} & \dots\dots & b_{-1}  }$$ 
tells us that $M/JM \cong S(a_1)$,
$JM/J^2M \cong S(a_0) \oplus S(b_0)^m$, $J^2M \cong S(b_{-1})^m$, and
$J^3 M = 0$.  In particular, this implies $\gamma _0  \alpha_{1
0}  = \varepsilon _0 \alpha _{1 i} = \varepsilon _{-1}\alpha _{0 i}\alpha
_{1 0} = \varepsilon _{-1} \beta _0 \alpha _{1 i} = 0$ for $1 \le i \le
m$.    Moreover, the graph communicates the existence of nonzero scalars
$k_1, \dots , k_m$ such that $\alpha _{0 i}\alpha _{1 0} = k_i \beta _0
\alpha _{1 i}$ for $1 \le i \le m$. 
(In Section 4, we will make the automatic choice $k_i = 1$, whenever we
encounter non-monomial relations of this type.)

We give a final example, this time of a non-projective $\la$-module $N$. 
That
$N$ has graph

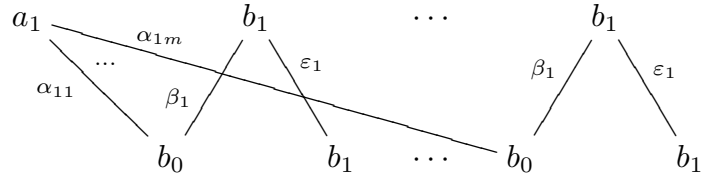
\ignore
$$\xymatrixcolsep{1pc}\xymatrixrowsep{2.5pc}
\xymatrix{ 
a_1 \edge[drr]_(0.4){\alpha _{1 1}}^(0.4){\dots}
\edge[drrrrrr]^(0.25){\alpha_{1 m}} &&& b_1
\edge[dl]_(0.65){\beta_1} \edge[dr]^(0.4){\varepsilon _1} &&
{\displaystyle{\cdots}} && b_1
\edge[dl]_(0.45){\beta _1} \edge[dr]^(0.45){\varepsilon _1} \\
 && b_0 && b_1 & {\displaystyle{\cdots}} & b_0 && b_1
 }$$ 
\endignore

\noindent relative to a sequence $x_0, x_1, \dots , x_m$ of top elements
(generating
$N$ and $K$-linearly independent modulo $JN$) means: $N/JN
\cong S(a_1) \oplus S(b_1)^m$ with $x_0$ generating a copy of $S(a_1)$
modulo $JN$ and each $x_i$ for $i \ge 1$
generating a copy of $S(b_1)$ modulo $JN$;  in other words,
$x_0$ is a top element of type $a_1$, while $x_1, \dots , x_m$ are top
elements of type $b_1$.  Moreover, we glean 
$JN/J^2N 
\cong S(b_0)^m \oplus S(b_1)^m$, showing 
$\alpha _{1 0}x_0 = 0$, and $J^2 N = 0$, which entails 
$$\beta _0 \beta _1 x_i = \varepsilon _0 \beta _1 x_i = \varepsilon _1^2
x_i = \beta _1 \varepsilon _1 x_i = 0 \quad \text{for} \,\, 1 \le i \le
m.$$  Finally, the graph holds the information that
$\alpha _{1 i}x_0 = k_i\beta _1 x_i$ for suitable scalars $k_i \in K^*$.

A final convention: A layered graph of a module $M$ is said to be a {\it
tree} if the underlying unlayered graph is a tree. The graph
of $N$ above is an example.

\head 3.  Stacks of algebras \endhead

We return to the situation where $\la = K\Q/I$ is an arbitrary finite
dimensional path algebra modulo relations.  For ease of notation, we
write $E$ for the set of primitive idempotents of $\la$ (= vertices of
$\Q$).  

\definition{Definition 2} 

\noindent {\rm (1)}  A {\it stacking partition\/} of $\Lambda$ is a
disjoint partition
$E = E' \cup E''$ satisfying the following two conditions: 

 (a) Every arrow of $\Q$ that starts in $E'$ also ends in $E'$, i\.e.,
$E'' \cdot (K\Q) 
\cdot E' = 0$. 

 (b)  Suppose $\alpha$ is an arrow in $\Q$ which starts in $E''$ and ends
in $E'$, and let $\beta$ be any arrow.  Then $\alpha \beta
\not\in I$ forces $\beta$ to start in a source of $\Q$.
\medskip

\noindent {\rm (2)}  Given a stacking partition of $\la$ as under (1), we
let $\be'$ (resp., $\be''$) be the sum of all idempotents in $E'$ (resp.,
$E''$), and set $\la' = \be' \la \be'$, $\la'' = \be''\la \be''$.

If both $\be'$ and $\be''$ are nonzero, we call  $\la$ a $2$-{\it
stack\/}, or say that $\la$ is obtained by {\it stacking
$\la''$ on top of\/} $\la'$. 
\medskip

\noindent {\rm (3)}  Suppose $\Delta_i$ for $0 \le i \le d$ are finite
dimensional algebras.  We call $\la$ a ($d+1$)-{\it stack of\/} $\Delta_0,
\dots, \Delta_d$ if there exist algebras $\la_0, \dots, \la_d$ with $\la_0
= \Delta_0$ and 
$\la_d = \la$, such that $\la_i$ is a
$2$-stack with $\Delta_i$ stacked on top of $\la_{i-1}$ for $i \ge 1$. 
\enddefinition

A more suggestive rendering of a 2-stack $\la$ obtained by
stacking $\la'' = K\Q''/I''$ on top of $\la' = K\Q'/I'$ is as follows: 
$\la = K\Q/I$, where $\Q$ is a quiver of the form 

\ignore
$$\xymatrixcolsep{1pc}\xymatrixrowsep{0.25pc}
\xymatrix{
{}\save[0,0]+(0,3);[0,2]+(0,-3)**\frm{-}\restore &\Q'' & \\
\ar[ddd]<1ex>_{\alpha_1} & \ar@{}[ddd]|{\displaystyle\cdots} &
\ar[ddd]<-1ex>^{\alpha_m} \\  \\  \\  
 && \\
{}\save[0,0]+(0,3);[0,2]+(0,-3)**\frm{-}\restore &\Q' & {}
 }$$
\endignore

\noindent Here $\alpha_1, \dots, \alpha_m$ are `new' arrows in $\Q$, and
$I$ is an admissible ideal containing the paths $\alpha_i \beta$ for all
arrows $\beta$ in
$\Q''$ starting in non-sources such that, moreover, $I \cap K\Q' = I'$
and
$I \cap K\Q'' = I''$.  A $3$-stack of algebras $\Delta_0$, $\Delta_1$,
$\Delta_2$ can thus be roughly visualized in the form

\ignore
$$\xymatrixcolsep{1pc}\xymatrixrowsep{0.25pc}
\xymatrix{
{}\save[0,0]+(0,3);[0,2]+(0,-3)**\frm{-}\restore &\Delta_2 & \\
\ar[ddd]<2ex> \ar[ddd]<3ex> & \ar@{}[ddd]|{\displaystyle\cdots} &
\ar[ddd]<-3ex> \ar[ddd]<-2ex> \\  \\  \\  
 && \\
{}\save[0,0]+(0,3);[0,2]+(0,-3)**\frm{-}\restore &\Delta_1 & \\
\ar[ddd]<2ex> \ar[ddd]<3ex> & \ar@{}[ddd]|{\displaystyle\cdots} &
\ar[ddd]<-3ex> \ar[ddd]<-2ex> \\  \\  \\  
 && \\
{}\save[0,0]+(0,3);[0,2]+(0,-3)**\frm{-}\restore &\Delta_0 & 
 }$$
\endignore

\noindent which motivates the terminology `stacking'.

In case
$\la$ is given in terms of a quiver and a generating set for $I$, the
problem of recognizing stacking partitions of $E$ can often be resolved
by mere inspection of the data. In general,
$E$ will have many different stacking partitions, a fact that can be used
to advantage in obtaining a maximum of homological information about
$\la$; this is witnessed by our applications.  To motivate condition (b)
in the definition of a stacking partition, we point to the following
obvious fact:  Whenever $M$ is a submodule of the radical of a projective
module, we have $e(M/JM) = 0$ for all sources $e$ of $\la$; in other
words, the set
$E(M)$ of Theorem 1 consists of non-sources in this situation.  

Moreover, we will see that condition (b) ensures that, on the level of
second syzygies, one obtains good separation of the `contributions' from
the layers of a
$2$-stack.  This condition aims specifically at applications involving
monomial algebras as building blocks, since their second syzygies become
structurally transparent.  If only the third or higher syzygies over the
algebras one wishes to stack are known to have good properties, it is
advantageous to relax our key definition as follows:  For $c \in
\NN$, a partition
$E = E' \cup E''$ is a {\it stacking partition of complexity $c$\/} in
case (a) holds and condition (b) is relaxed as follows:  Given any arrow
$\alpha \in  E' \cdot (K\Q) \cdot E''$, the product $\alpha \beta$ belongs
to
$I$ for all arrows $\beta$ starting in  the endpoint of a path of length
$c$ in $\Q$. In this sense, our stacking partitions are of complexity
$1$.      

For the remainder of this section, we assume that $E = E' \cup E''$ is a
stacking partition of $\la$ with both $E'$ and $E''$ nonempty. We retain
the notation $\be'$, $\be''$, $\Lambda '$, and
$\Lambda''$ from part (2) of the definition.   Moreover, we let $\Q'$ and
$\Q''$ be the full subquivers of $\Q$ with vertex sets $E'$ and
$E''$, respectively.   Clearly, 
$$\la' \cong K\Q'/ I\cap K\Q' \ \ \ \ \
\text{and} \ \ \ \ \  \la'' \cong K\Q''/ I\cap K\Q'',$$ 
and the Jacobson
radicals of these algebras are $J' = \be'J\be'$ and
$J'' = \be''J\be''$, respectively.  On the side, we note that, whenever
$\Q''$ contains a loop, $\la$ does not result from iterated one-point
extensions of $\la'$.

Observe that, for any $N \in \Lamod$, the $\la'$-component $\be' N$ is a
$\la$-submodule of $N$.  The analogous statement for $\be'' N$ is
obviously false; it even fails for first syzygies of $\la$-modules.  In
general, we only have a $K$-vectorspace decomposition $\Omega^1(N) = \be'
\Omega^1(N) \oplus \be'' \Omega^1(N)$.  However, on the level of second
syzygies in
$\Lamod$, we obtain nice splittings into $\la'$ and $\la''$-components. 
This is the pivotal point in relating the homological properties of stacks
to those of their building blocks.

\proclaim{Proposition 3}  If $X$ is a second syzygy of a left
$\la$-module $N$, then both $\be' X$ and $\be'' X$ are $\la$-submodules of
$X$, and thus the $K$-vector space decomposition 
$$X = \be'X \oplus \be''X$$ is a $\la$-direct sum. 

In particular, $\pdim_{\la} N < \infty \iff \pdim_{\la} \be'X < \infty
\text{\ and \ }\pdim_{\la} \be'' X < \infty$.

\endproclaim

\demo{Proof}  We only need to show that $\be''X$ is a $\la$-submodule of
$X$.  By hypothesis,
$X$ is the kernel of a map $f: P
\rightarrow M$, where $M$ is a submodule of the radical of a projective
module.  This entails $e(M/JM) = e(P/JP) = 0$ for all sources $e$ of $\Q$.
Consequently, condition (b) of a stacking partition implies the
following:  Whenever $\alpha$ is an arrow in $\Q$ and $x \in \be''X
\subseteq P$, we have
$\alpha x = \be'' \alpha x$.  Indeed, 
$x$ is a linear combination of elements $p_i x_i$, where the $x_i$ are top
elements of $P$ and the
$p_i$ paths of positive lengths starting in non-sources and ending in
$E''$; hence each product $\alpha p_i$ is either zero in $\la$ or
else a path ending in $E''$. \qed \enddemo

This focuses the discussion on the question of how the projective
dimensions of the components $\be' X$ and $\be'' X$, viewed as $\la'$- and
$\la''$-modules respectively, relate to their $\la$-projective
dimensions.  For $\be' X$ this is obvious  --  we will nevertheless record
it  --  for
$\be'' X$ it is a far more intricate problem. 

The following straight-forward lemma only uses the fact that $\be'' \la
\be' = 0$.

\proclaim{Lemma 4}  

\noindent {\rm(1)} If $N$ is any left $\la$-module, the minimal projective
resolution of $\be' N$ over $\la'$ coincides with the minimal projective
resolution of $\be' N$ over $\la$, i\.e\., $\Omega^i_\la (\be' N) =
\Omega^i_{\la'} (\be' N)$ for all $i$.  (Note, however, that, given a
projective
$\la$-module $Q$, its
$\la'$-component
$\be' Q$ need not be projective as a $\la$- or, equivalently, as a
$\la'$-module.)

\noindent {\rm(2)}  Given any left $\la$-module $N$ and a
$\la$-projective cover $g: Q \rightarrow N$ with kernel $M$, the
restriction $\be''g: \be''Q \rightarrow \be''N$ is a $\la''$-projective
cover of $\be''N$ with kernel $\be''M$.  

In particular, $\Omega^i_{\la''}(\be''N) = \be'' \Omega^i_{\la} (N)$ for
all
$i \ge 0$.
\endproclaim

\demo{Proof}  The first claim is obvious.  For part (2), we observe that
the algebra $\la''$, viewed as a left module over itself, has the
following decomposition into left ideals, by the definition of a stacking
partition:  Namely
$\la'' =
\bigoplus_{e \in E''} e \la \be'' = \bigoplus_{e \in E''} e \la =
\bigoplus_{e \in E'',\ f \in E} e \la f$.  That $\be'' J = \be''J \be''$
and hence
$\be''M \subseteq J''(\be''Q)$, has similar reasons.  \qed
\enddemo

Noting that Proposition 3 carries over to direct summands of second
syzygies in $\Lamod$, we derive the following consequence.

\proclaim{Proposition 5}

\noindent {\rm(1)}  For any left $\la$-module $N$, 
$$\pdim_{\la} \be' N = \pdim_{\la'} \be'N.$$

\noindent {\rm(2)}  Suppose that $X$ is any $\la$-direct summand of a
second syzygy in $\Lamod$. Then
\smallskip

{\rm (a)} $\pdim_{\la} \be'' X \ge \pdim_{\la''} \be'' X$.
\smallskip

{\rm (b)} If $\pdim_{\la} \be'' X < \infty$, then
$$\pdim_{\la} \be''X \le \pdim_{\la''} \be''X + \bold{t} + 1,$$ where
$\bold{t} = 
\max\{\pdim_{\la'} \be'\la e \mid e \in E''\setminus\{\text{sources of }
\Q'' \} \text{\ and\ }   \pdim_{\la'} \be'\la e < \infty\}$ in case the
relevant set is nonempty, and $\bold{t} = -1$ otherwise.
\smallskip

{\rm (c)}  If $\pdim_{\la} \be'' X = \infty$, then either $\pdim_{\la''}
\be'' X = \infty$, or else there exists an idempotent $e \in E''$ such
that $\pdim_{\la'} \be'\la e =\infty$. 
\endproclaim

\demo{Proof}  (2)  From Proposition 3 we know that $\be''X$ is a
$\la$-submodule of $X$.  Thus (a) is an immediate consequence of Lemma 4.

For part (b), suppose $\pdim_{\la} \be''X < \infty$.  It is, moreover,
harmless to assume that $\pdim_{\la''} \be''X < \infty$.  If $f_0: P_0
\rightarrow \be''X$ is a $\la$-projective cover, then clearly
$e(P_0/JP_0) = 0$ for all sources $e$ of $\Q$, and hence the simple
summands of $P_0/JP_0$ correspond to idempotents in  
$E''\setminus\{\text{sources of }
\Q'' \}$.  Using once more Lemma 4 and Proposition 3, we further see that
$\Omega^1_{\la''}(\be''X) = \be''
\Omega^1_{\la}(\be'' X)$ is a $\la$-submodule of $\Omega^1_{\la}(\be''
X)$ and, in view of $\be'P_0 \subseteq \ker(f_0)$, we obtain a
$\la$-direct decomposition 
$$\Omega^1_{\la}(\be'' X) = \be' P_0 \oplus \Omega^1_{\la''}(\be''X).$$
Our hypothesis that $\pdim_{\la} \be'' X$ be finite ensures that
$\pdim_{\la'}
\be' P_0 \le \bold{t}$.  So, if
$\pdim_{\la''} \be'' X \allowmathbreak = 0$, the desired inequality
follows.  Otherwise, we repeat the preceding argument with
$\Omega^1_{\la}(\be'' X)$ instead of
$X$, to obtain
$$\Omega^2_{\la}(\be'' X) = \be' P_1 \oplus \Omega^2_{\la''}(\be''X),$$
where $P_1$ is a $\la$-projective cover of $\be''\Omega^1_{\la}(\be''
X)$.  (This is a legitimate move, because $\be'' X$ is a $\la$-direct
summand of
$X$, and hence $\Omega^1_{\la}(\be'' X)$ again satisfies the blanket
hypothesis of (2).)  Thus our claim is also true in case
$\pdim_{\la''}\be'' X
\allowmathbreak = 1$.  An obvious induction now completes the argument.

Part (c) is obtained analogously. \qed\enddemo

It is now easy to deduce that the little finitistic dimensions of stacks
are governed by the inequalities mentioned in the introduction. 
Since
$$\la = \pmatrix \la'' & 0 \\ \be' \la \be'' & \la' \endpmatrix,$$ 
these
inequalities actually hold in far greater generality, as 
the methods of \cite{\Fos, Corollary 4.21} show.  For the reader's
convenience, we include an argument for our special case. 

\proclaim{Corollary 6} The left little finitistic dimension of $\la$
satisfies the inequalities  
$$\findim \la'  \le \findim \la \le  \findim \la' +  \findim
\la'' + 1.$$   Moreover, 
$$\findim_n \la'  \le  \findim_n \la$$  for all $n \in \NN$.   
\endproclaim

\demo{Proof} The final set of inequalities is obvious, as is the fact
that $\findim \la' \le \findim \la$.  To check the upper bound on
$\findim \la$, let $N \in \lamod$ with $\pdim_{\la} N <
\infty$.  Repeated use of Lemma 4 yields $\Omega^2_{\la''}(\be'' N)$
$=$ $\Omega^1_{\la''}\bigl( \be'' \Omega^1_{\la} (N) \bigr)$ $=$ $\be''
\Omega^2_{\la}(N)$ $=$ $\be'' X$ if $X = \Omega^2_{\la}(N)$.  By
induction, 
$$\Omega^k_{\la''}(\be'' N) = \be'' \Omega^k_{\la}(N) = \be''
\Omega^{k-2}_{\la}(X)$$ for all $k \ge 2$.  In particular, $\be'' N$ has
finite projective dimension over $\la''$, which shows $\pdim_{\la''}
\be'' N \le
\findim \la''$.  We first deal with the case where $\findim \la'' = 0$. 
In this situation, $\pdim_{\la''} (\be''N) = 0$, whence $\be''
\Omega^1_{\la} (N) = 0$ by Lemma 4; this makes $\Omega^1_{\la}(N)$ a
$\la'$-module with
$\pdim_{\la}\Omega^1_{\la}(N) = \pdim_{\la'} \Omega^1_{\la}(N)$ and thus
ensures $\findim \la \le \findim \la' + 1$.   

Now suppose that $\findim \la'' \ge 1$.  Without loss of generality, we
may assume that
$\pdim_{\la''} \be'' N \ge 2$, for otherwise $\be''X$ would vanish and
$\pdim_{\la}N$ would be bounded above by the sum $2+ \pdim_{\la}
\be' X$ $=$ $2 + \pdim_{\la'} \be' X$, which is at most $\findim \la' +
2$.  Hence 
$\pdim_{\la} N = 2 + \max\{\pdim_{\la} \be' X, \pdim_{\la} \be'' X\}$.
Now $\pdim_{\la} \be' X = \pdim_{\la'} \be' X$ is bounded above by
$\findim
\la'$, and
$\pdim_{\la} \be'' X \le \pdim_{\la''} \be'' X + \bold{t} + 1$ by
Proposition 5.  The final term equals $\bigl(\pdim_{\la''} \be''N - 2
\bigr) +
\bold{t} + 1$ which is in turn bounded from above by $\bigl(\findim
\la'' - 2 \bigr) + \findim \la' + 1$.  Consequently,
$$\pdim_{\la}N \le \findim \la' + \findim \la'' + 1$$ as required.
 \qed \enddemo

In general, $\findim \la''$ is not a lower bound for $\findim \la$  -- 
examples to the contrary are ubiquitous.

\example{Example 7}  Suppose that $\la$ is a monomial algebra such that
the indecomposable projectives in
$\lamod$ have the following graphs:
\medskip
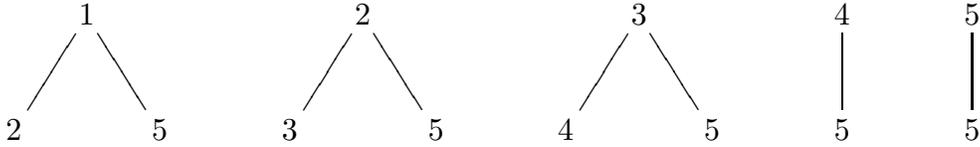
\ignore
$$\xymatrixcolsep{1pc}\xymatrixrowsep{2pc}
\xymatrix{
 &1 \edge[dl]\edge[dr] &&&&2 \edge[dl]\edge[dr] &&&&3
\edge[dl]\edge[dr] &&&4 \edge[d] &&5 \edge[d] \\ 2 &&5 &&3 &&5 &&4 &&5
&&5 &&5  }$$
\endignore
\medskip

Then the vertex sets $E' = \{5\}$ and $E'' = E \setminus E'$ define a
stacking partition of $\la$,  with $\findim \la = \findim
\la' = 0$, and $\gldim \la'' = 3$. \endexample  

Already the case where the underlying graph of $\Q$ is a  Dynkin diagram
of type $A_n$ yields examples showing the upper bound on $\findim \la$ of
Corollary 6 to be optimal, even when all algebras involved have finite
global dimensions.  

\example{Example 8}  Suppose that $\Q$ has underlying graph $A_5$ with
arrows $i \rightarrow i+1$, and set $\la = K\Q/I$, where I is generated by
all paths of length $2$.   Consider the stacking partition $E' = \{4,5\}$
and $E'' = \{1,2,3\}$ of $\la$, and observe, that $\gldim \la' = 1$,
$\gldim \la'' = 2$, and $\gldim \la = \gldim \la' + \gldim \la'' + 1$.
\endexample

We add a few specialized comments addressing stacks of monomial algebras
for use in the next section.  Suppose that
$\la'' = K\Q''/(K\Q'' \cap I)$ is a monomial algebra, and, as in Theorem
1, consider the set of critical paths of $\la''$:
$$\multline 
\eusm{P''} = \{ p \in \Lambda'' \mid p \text{ is a nontrivial path of
positive length,} \\ 
\text{starting in a non-source of } \Q'', \text{ with\ }
\pdim_{\Lambda''}\Lambda'' p < \infty \}.
\endmultline$$    
Moreover, let $\bold{s''}$ be the supremum
of the projective dimensions attained on this set, as in Theorem 1. 
Since, given any $N \in \Lamod$, the $\la''$-module $\be'' N$ has second
syzygy $\Omega^2_{\la''}(\be'' N) \cong \be'' \Omega^2_{\la}(N)$, that
theorem shows $\be'' \Omega^2_{\la}(N)$ to split into a direct sum of
cyclic modules isomorphic to left ideals $\la'' q$  for suitable paths $q$
in $\la''$.  Clearly, the $\la'' q$ are even left ideals of $\la$, i\.e\.
$\la''q = \la q$, since the eligible paths $q$ start in non-sources. 
Yet, in general, the
$\la$-projective dimensions of the $\la'' q$ will still exceed their
$\la''$-projective dimensions; in particular, the former may be infinite
while the latter are finite (see Example 7 above).  In the construction
we are targeting, this problem does not arise, however, since all of
the end points of the paths in
$\eusm{P''}$ are  `homogeneous' in the following sense:  A vertex $e \in
E''$ is called {\it homogeneous\/} in case all arrows of $\Q$ starting in
$e$ end in
$E''$.  If all paths in $\eusm{P''}$ end in homogeneous vertices, then
clearly  $\pdim_{\la} \la''p = \pdim_{\la''} \la''p$ for all $p \in
\eusm{P''}$.  Combining these considerations with the preceding results,
we obtain:

\proclaim{Corollary 9}  Let $\la''$ be a monomial algebra.  Retaining the
above notation, suppose that all end points of the paths in
$\eusm{P''}$ are homogeneous.   If $N$ is any left $\la$-module of finite
projective dimension, then the $\la''$-projective dimension of the syzygy
$\Omega^2_{\la''} (\be'' N)$ coincides with its $\la$-projective
dimension.

In particular:  If $N$ is finitely generated and $\pdim_\la N > \findim
\la' + 2$, then $\pdim_{\la''}
\be''N \allowmathbreak = \pdim_\la N$. 
 \endproclaim

\demo{Proof}  Proposition 3 yields a $\la$-direct decomposition 
$$\Omega^2_\la (N) = \be' \Omega^2_\la(N) \oplus \be'' \Omega^2_\la (N),$$
which shows $\Omega^2_{\la''} (\be''N) = \be'' \Omega^2_\la (N)$  to have
finite projective dimension as a
$\la$-module.  By Proposition 5, we infer that $\Omega^2_{\la''} (\be''N)$
has finite projective dimension also as a $\la''$-module, and
consequently Theorem 1 shows this syzygy to decompose in the form 
$$\Omega^2_{\la''} (\be''N) = \bigoplus_{i \in I} \la'' p_i,$$
where the $p_i$ are paths in $\eusm P''$.  Using our homogeneity
hypothesis, we thus conclude $\pdim_{\la''} \Omega^2_{\la''} (\be''N) =
\pdim_\la \Omega^2_{\la''} (\be''N)$ as required.  For the final
assertion, one just has to keep in mind that the $\la$-projective
dimension of $\be' \Omega^2_\la (N)$ equals the $\la'$-projective
dimension (Proposition 5(1)).  \qed
\enddemo

\head 4. The Key Examples \endhead

\proclaim{Theorem 10}  Given any  increasing function
$f:\Bbb{N}\rightarrow\NN_{\ge 2}$ that takes only finitely many
distinct values, there exists a connected finite dimensional algebra $\la
= K\Q /I$ with vanishing radical cube such that $$\roman{fin} \,\,
\roman{dim}_n\la = f(n)\,\, \roman{for}\,\, n \in \Bbb{N}.$$ Moreover, if
$d = \max f - \min f,$ then $\la$ can be constructed as a $\left( \lceil
d/2 \rceil +1 \right)$-stack of monomial algebras with the additional
property that, for each $n \in \Bbb{N}, \,\, \roman{fin} \,\,
\roman{dim}_n\la$ is attained on an $n$-generated module of Loewy length
at most $2$ having a tree graph.
\endproclaim

Our proof consists of a recursive construction technique involving $d$
stages, where $d$ is as in the theorem.  To avoid overly cumbersome
notation, we will only deal with two types of step functions, the first
exhibiting a single jump of size $s,$ the other involving two jumps of
sizes $s$ and $t,$ respectively.  The general recursive pattern is clear
from our constructions.

\subhead Part I of the proof of Theorem 10 \endsubhead 
 Fix $m, \, r, \,
s \in \Bbb{N}$ with $m,r \ge 2,$ and assume that $f(k) = r$ for $k \le m
- 1$, while $f(k) = r+s$ for $k \ge m.$

We start by constructing a sequence $\la _0, \la _1, \dots , \la _s$ of
finite dimensional algebras, where $\la _0$ is a monomial algebra and
each $\la _\ell$ for $\ell \ge 1$ results from stacking another monomial
algebra on top of $\la _{\ell - 1}.$  Then we will show that, for $0 \le
\ell \le s$,
$$\findim_k \la _\ell = r \ \ \text{for} \ \ k \le m - 1 \quad \text{and}
\quad  \findim_k \la _\ell = r + \ell \ \ \text{for} \ \ k \ge m.
\tag{$\dagger$} $$

Note that, for the present choice of $f$, the difference $\max f - \min
f$ equals $s.$  When we assemble the building blocks for our argument
under the `final claim of Part I', we will explain how the algebra $\la =
\la _s$ can be obtained through a more economical stacking of only
$\lceil s/2
\rceil + 1$ monomial layers.

\remark{Step $0$}  We base $\la_0$ on the following quiver $\Q
^{(0)}$:

$$
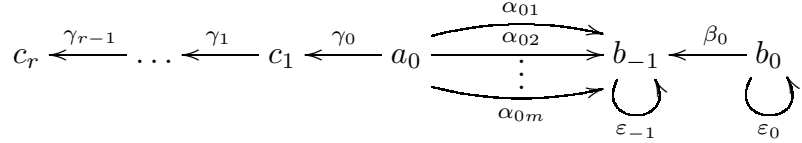
\xymatrix{
  c_r & \dots \ar[l]_{\gamma _{r-1}} & c_1 \ar[l]_{\gamma _1} & a_0
\ar[l]_{\gamma _0} \ar@/^/[rr]<0.85ex>^{\alpha_{01}}
\ar[rr]^{\alpha_{02}} \ar@/_/[rr]<-1.8ex>^{\vdots}_{\alpha _{0m}} & &
b_{-1} \ar@(dl,dr)_{\varepsilon _{-1}} & b_0 \ar[l]_(0.4){\beta _0}
\ar@(dl,dr)_{\varepsilon _0} }$$

\noindent The algebra $\la _0$ results from $K\Q ^{(0)}$ by factoring out
the ideal generated by the relations $\varepsilon _i ^2$ for $i = -1,
0$, as well as $\varepsilon _{-1}\beta_0$, $\beta _0\varepsilon _0$ and
$\gamma _j
\gamma _{j-1}$ for $1 \le j \le r - 1$.  Thus the indecomposable
projective left $\la _0$-modules $\la _0 a_0, \,\, \la _0 b_i,$ and $\la
_0 c_i$ have graphs

$$\xymatrixcolsep{0.75pc}
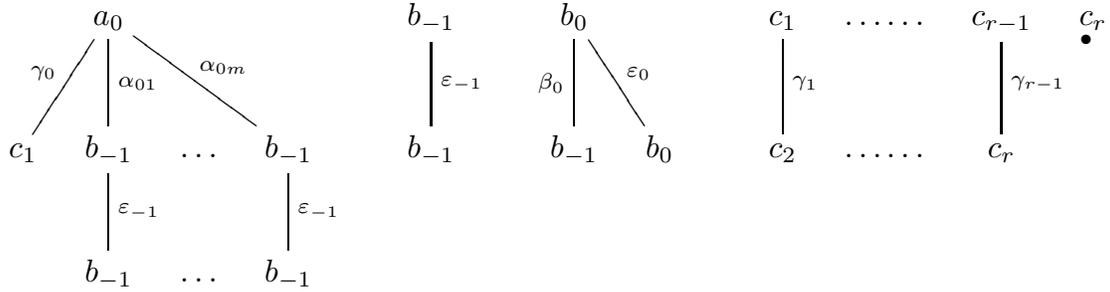
\xymatrix{ & a_0 \edge[dl]_{\gamma _0} \edge[d]^{\alpha _{01}}
\edge[drr]^{\alpha _{0m}} & & & & b_{-1} \edge[d]^{\varepsilon
_{-1}} & & b_0 \edge[d]_{\beta _0} \edge[dr]^{\varepsilon_0} & & & c_1
\edge[d]^{\gamma _1} & \dots\dots & c_{r-1} \edge[d]^{\gamma _{r-1}} &
\underset {\bullet}\to c_r \\ c_1 & b_{-1} \edge[d]^{\varepsilon _{-1}} &
\dots & b_{-1} \edge[d]^{\varepsilon _{-1}} & & b_{-1} & & b_{-1} & b_0 &
& c_2 & \dots\dots & c_r \\ & b_{-1}& \dots & b_{-1}  }$$

\noindent respectively.   

\endremark

\remark{Step \rom{$\ell$} \rom{(for $\ell \ge 1$)}}  The algebra $\la
_\ell$ is based on a quiver $\Q ^{(\ell)}$ having vertex set $\Q
_0^{(\ell)} = \Q _0^{(\ell - 1)}\cup \{ a_{\ell}, b_{\ell} \}$ and
additional arrows 
$$\alpha _{\ell 0}: a_{\ell} \rightarrow a _{\ell - 1}, \quad \alpha
_{\ell 1}, \dots , \alpha _{\ell m}:a_{\ell}\rightarrow b_{\ell - 1},
\quad \beta _{\ell}:b_\ell \rightarrow b_{\ell - 1}, \quad \varepsilon
_\ell :b_\ell \rightarrow b_\ell.$$

The graphs of the indecomposable projective $\la _\ell$-modules $\la
_\ell e$, where $e$ is a vertex of $\Q ^{(\ell - 1)}$, are the same as
those of the corresponding $\la_{\ell-1}$-modules $\la _{\ell - 1}e$, and
the `new' indecomposable projective
$\la _\ell$-modules $\la _\ell a_\ell$ and $\la _\ell b_\ell$  have
graphs

$$\xymatrixcolsep{1.5pc}
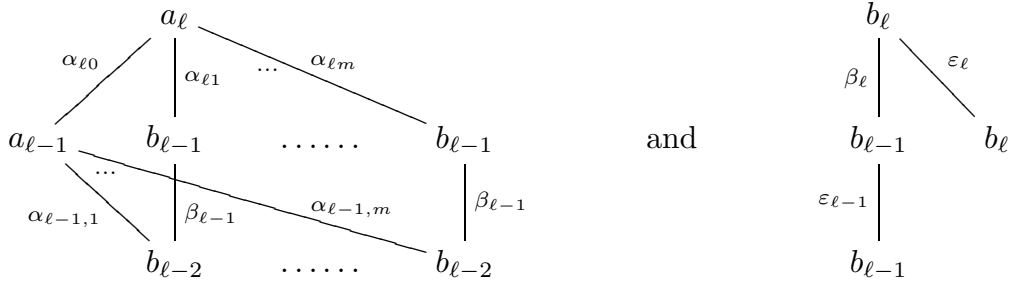
\xymatrix{ & a_\ell \edge[dl]_{\alpha_{\ell 0}} \edge[d]^{\alpha _{\ell
1}} \edge[drr]^{\alpha _{\ell m}}_(0.35){\dots} & & & & & & b_\ell
\edge[d]_{\beta _\ell} \edge[dr]^{\varepsilon _\ell} \\ a_{\ell - 1}
\edge[dr]_{\alpha_{\ell -1, 1}}^(0.35){\dots} \edge[drrr]^(0.7){\alpha
_{\ell -1,m}} & b_{\ell - 1} \edge[d]^(0.6){\beta _{\ell -1}} &
\dots\dots & b_{\ell - 1} \edge[d]^{\beta _{\ell - 1}} & & \text{and} & &
b_{\ell -1} \edge[d]_{\varepsilon _{\ell - 1}} & b_\ell \\ & b_{\ell - 2}
& \dots\dots & b_{\ell - 2} & & & & b_{\ell - 1} }$$

\noindent respectively; here each subgraph of the form
$$\xymatrixcolsep{0.6pc} \xymatrixrowsep{0.6pc}
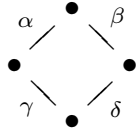
\xymatrix{ & \bullet \edge[dl]_{\alpha} \edge[dr]^{\beta} \\
\bullet \edge[dr]_{\gamma} & & \bullet \edge[dl]^{\delta} \\ & \bullet 
}$$
\newline of the graph of $\la _\ell a_\ell$ is to reflect a relation
$\gamma \alpha - \delta \beta$ of $\la _\ell$.  (In Section 2, we
explained how the graphs of $\la _1 a_0$ and $\la _1 a_1$ are to be
interpreted.)

Our proof for the asserted values of the finitistic dimensions of the
$\la _\ell$ (see ($\dagger$) above) will be deduced from the following
sequence of lemmas.  As a prerequisite for applying the results of
Section 3, we only need the obvious fact that each $\la_\ell$ is an
$\ell$-stack of monomial algebras; if $E = \Q _0^{(\ell)}$ is the vertex
set of $\Q ^{(\ell)}$, we use the stacking partition $E' = \Q
_0^{(\ell - 1)}$ and $E'' = \{ a_\ell , b_\ell \}$ of $E$.  In particular,
this means that
$\la _k$-mod embeds into $\la _{\ell}$-mod for $k <  \ell$; that, for any
$\la _k$-module
$X$, the $\la _\ell$-projective cover of $X$ coincides with the $\la
_k$-projective cover (Lemma 4), and hence that $\pdim_{\la _k}X =
\pdim_{\la _\ell}X$.  Consequently, there is no need to emphasize the base
algebra in dealing with projective covers and computing projective
dimensions.  We will denote the Jacobson radical of $\la _\ell$ by
$J_\ell$.  Note that
$J_k X = J_\ell X$ for all $\la _k$-modules $X$ by Section 3.

\proclaim{Lemma 11}  Let $\ell \ge 0$, and suppose that $N$ is a $\la
_{\ell +1}$-module of finite projective dimension.  Then $b_{\ell
+1}\Omega^1(N) = 0$; in particular, $\Omega^1(N)$ is a $\la _\ell$-module
and $\findim\la _{\ell + 1} \le 1 + \findim\la _{\ell}$.
\endproclaim

\demo{Proof} If Q is a projective cover of $N,$ write $Q = Q_1 \oplus
Q_2$, where $Q_2$ is a direct sum of copies of $\la _{\ell +1}b_{\ell
+1}$ and $Q_1$ has no direct summands in common with $Q_2.$  Then
$b_{\ell +1}JQ = b_{\ell +1}JQ_2 = \varepsilon _{\ell +1}Q_2$ is a direct
sum of copies of $S(b_{\ell +1})$ and a $\la _{\ell +1}$-direct summand
of $JQ$.  Since $\pdim S(b_{\ell +1})=\infty$, this implies $b_{\ell
+1}\Omega ^1(N) = 0$.  Using $\alpha_{\ell + 1}J_{\ell + 1} = 0$, we infer
that $\Omega ^1(N)$ is indeed a $\la _\ell$-module.  The final claim is
thus obvious. \qed   
\enddemo

\proclaim{Lemma 12}  Let $\ell \ge 0$, and suppose $N$ is a $\la _{\ell
+1}$-module of finite projective dimension containing a top element $y$
of type $a_{\ell +1}$ such that $\alpha _{\ell +1, 0} y = 0$.  Then
$N/JN$ contains the simple module $S(b_{\ell +1})$ with multiplicity at
least $m$.  
\endproclaim

\demo{Proof}  Let $\psi :Q = \la_{\ell +1}\tilde y \oplus \tilde Q
\rightarrow N$ be a projective cover, where $\tilde y$ is a top element
of type $a_{\ell +1}$ of $Q$, such that $\psi (\tilde y) = y$.  Moreover,
set $M = \ker (\psi )$, i.e., $M \cong \Omega ^1(N)$.  Then $M$ is a
$\la _\ell$-module by Lemma 11, and $x = \alpha _{\ell +1, 0}\tilde y$ is
a top element of type $a_\ell$ of $M$.  Finally, let $\phi : P = \la
\tilde x
\oplus \tilde P \rightarrow M$ be a projective cover of $M$, where
$\tilde x \in P$ is a top element of type $a_{\ell}$ with $\phi (\tilde
x) = x$.

We focus on the case $\ell \ge 1$ in the sequel, the argument for $\ell =
0$ being analogous, modulo small adjustments. (Due to the difference in
make-up of $\la _\ell a_\ell$ for $\ell \ge 1$ and $\ell = 0$, in the
latter case, the element $\varepsilon _{-1}\alpha_{0i}\tilde x$ takes
over the role played by $\beta _{\ell - 1}\alpha _{\ell i} \tilde x$
below.)

Fix $i \in \{ 1, \dots ,m\}$.  Since $x = a_\ell x$ belongs to $J_{\ell
+1}\tilde y \cong J_{\ell +1} a_{\ell +1}$, we see that $\alpha _{\ell
i}x$ is nonzero and generates a copy of $S(b_{\ell - 1})$ in the socle of
$M$, whereas $\beta _{\ell - 1}\alpha_{\ell i}x$ is zero.  Consequently,
$\alpha _{\ell i}\tilde x \notin \ker (\phi ) = \Omega ^2(N)$, whereas
$\beta _{\ell -1}\alpha _{\ell i}\tilde x \in \soc \ker (\phi )\setminus
\{0\}$.  If $\beta _{\ell - 1}\alpha _{\ell i}\tilde x$ were a top
element of $\ker (\phi )$, necessarily of type $b_{\ell - 2}$, this
element would generate a direct summand isomorphic to $S(b_{\ell - 2})$
in $\ker (\phi ) = \Omega ^2(N)$, clearly an impossibility.  Therefore,
$\beta_{\ell - 1}\alpha _{\ell i}\tilde x$ lies in $J _{\ell}\ker
(\phi ) = J _{\ell - 1}\ker (\phi )$; in fact, $\beta_{\ell - 1}\alpha
_{\ell i}\tilde x$ belongs to $J _{\ell - 1}\ker (\phi ) \cap \la _{\ell
-1}\alpha_{\ell i}\tilde x$ and therefore also to
$J_{\ell - 1}\pi (\ker (\phi )) \subseteq J_{\ell -1}J_{\ell}\tilde x$,
where $\pi : P \rightarrow \la \tilde x$ is the projection along $\tilde
P$.  Inspection of the graphs of the indecomposable projective $\la
_\ell$-modules thus reveals the existence of a top element $z_i$ of
$J_{\ell -1}P$   with the following properties: (1) $\alpha _{\ell i}
\tilde x - z_i$ is a top element of $\ker (\phi)$ of type $b_{\ell - 1}$,
and (2) $\beta_{\ell - 1}(\alpha_{\ell i}\tilde x - z_i) = \beta _{\ell -
1}\alpha _{\ell i}\tilde x$.  By (1), the top element $z_i$ of $J_{\ell -
1}P$ is of type $b_{\ell - 1}$.  From (2) we infer that $z_i = \beta
_\ell \tilde x_i$ for some top element $\tilde x_i$ of type $b_\ell$ of
$P$.  Thus $x$ and $x_i = \phi (\tilde x_i)$ are top elements of $M$
having types $a_\ell$ and $b_\ell$, respectively, with
$\alpha_{\ell i}x = \beta _\ell x_i$.  In particular, $\beta _\ell x_i$
belongs to $M \cap \la _{\ell +1}\tilde y \subseteq \la _{\ell +1}a_{\ell
+1}$.  So, if $\sigma :Q \rightarrow \la _{\ell +1}\tilde y$ is the
projection along $\tilde Q$, then $\alpha _{\ell i} x = \beta _\ell
\sigma (x_i)$, and the structure of $\la _{\ell + 1}a_{\ell +1}$
guarantees $\sigma (x_i) = \alpha _{\ell +1, i}\tilde y$.  From $\beta
_\ell (x_i - \sigma (x_i)) = 0$ we moreover glean $x_i - \sigma (x_i) =
\beta _{\ell + 1} \tilde y_i$, where $\tilde y_i$ is either zero or a top
element of type $b_{\ell + 1}$ of $Q$.

Next we check that $\tilde y_1, \dots, \tilde y_m$ are actually top
elements of $Q$ which are $K$-linearly independent modulo $JQ$: Indeed,
the $x_i = \sigma (x_i) + \beta _{\ell +1}\tilde y_i$ for $1 \le i \le m$
are linearly independent top elements of type $b_\ell$ of $M$, since the
multiples $\beta _\ell x_i$ are linearly independent by construction.  We
infer that the elements $$\varepsilon _\ell x_i = \varepsilon _\ell \beta
_{\ell +1}\tilde y_i, \,\, 1 \le i \le m$$ 
generate $m$ independent
copies of $S(b_\ell)$ in the socle of $M$; for otherwise we would obtain
$b_\ell \Omega ^1(M) \ne 0$, contradicting Lemma 11.  This forces $\tilde
y_1, \dots , \tilde y_m$ to be linearly independent modulo $JQ$ and thus
gives rise to $m$ top elements $y_i = \psi (\tilde y_i)$ of $N$ which are
linearly independent modulo $JN$. \qed
\enddemo

\endremark

\proclaim{Lemma 13} Let $\ell \ge 1$, and suppose $M$ is a $\la
_\ell$-module of Loewy length $2$ with $a_\ell M = 0$, but $b_\ell M \ne
0$.  Then $\roman{p} \,\, \roman{dim}\,\, M = \infty$.
\endproclaim

\demo{Proof} By hypothesis, the projective cover of $M$ does not contain
a summand $\la _\ell a _\ell$, but does contain a copy of $\la _\ell
b_\ell$.  Inspection of the graphs of 
$\la _\ell b_\ell$ and $\la _\ell b_{\ell - 1}$ thus makes it clear that
$\Omega ^1(M)$ either has a direct summand isomorphic to $S(b_{\ell -
1})$, or else a direct summand isomorphic to the module $X_{\ell - 1}$
with graph
$\vcenter{

\xymatrixrowsep{0.85pc} \xymatrixcolsep{1pc}
\xymatrix{ b_{\ell - 1} \edge[d] \\ b_{\ell - 1}  }}$. We know $S(b_{\ell
- 1})$ to have infinite projective dimension and compute $\Omega ^{\ell -
1}(X_{\ell - 1}) = X_0$.  Thus $\Omega ^\ell(X_{\ell - 1}) \cong
S(b_{-1})$, and we conclude that $X_{\ell - 1}$ has infinite projective
dimension as well.  \qed
\enddemo

\proclaim{Lemma 14} Let $\ell \ge 1$, and suppose that $N$ is an
indecomposable non-projective $\la _{\ell + 1}$-module of finite
projective dimension such that
$a_{\ell + 1}N \ne 0$ or $b_{\ell + 1}N \ne 0$.  Then $a_\ell \Omega
^1(N) \ne 0$.
\endproclaim

\demo{Proof}  By Lemma 13, it suffices to show that either $a_\ell \Omega
^1(N) \ne 0$ or $b_\ell \Omega ^1(N) \ne 0$.  Let $\psi : Q
\rightarrow N$ be a projective cover of
$N$, say $Q = Q_1 \oplus Q_2 \oplus Q_3$ where $Q_1$ is a direct sum of
copies of $\la _{\ell + 1}a_{\ell + 1}, \,\, Q_2$ a direct sum of copies
of $\la _{\ell + 1}b_{\ell + 1}$, and $Q_3$ has no direct summands in
common with
$Q_1 \oplus Q_2$.  From Lemma 11, we know $b_{\ell + 1}\Omega ^1(N) = 0$. 
Assume, to the contrary of our claim, that $a_\ell \Omega ^1(N) = b_\ell
\Omega ^1(N) = 0$, and note that $b_\ell \Omega ^1(N) = b_{\ell +1}
\Omega ^1(N) = 0$ implies $\Omega ^1(N) = \ker (\psi ) \subseteq Q_1
\oplus Q_3$, which places a direct summand isomorphic to $Q_2$ into $N$. 
In view of the fact that $N$ is indecomposable nonprojective, this
entails $Q_2 = 0$.  Moreover, our assumption forces $\ker (\psi )$ to be
contained in $b_{\ell - 1}Q_1 \oplus JQ_3$; however, $\ker (\psi )
\nsubseteq JQ_3$, due to indecomposability of $N$ - indeed $Q_2 = 0$
implies $Q_1 \ne 0$ by hypothesis.  Observe moreover that $b_{\ell -
1}Q_1 = \beta _\ell JQ_1$ is a nonzero direct sum of copies of $S(b_{\ell
- 1})$.  Again we deduce that $\ker (\psi )$ either contains a copy of
$S(b_{\ell - 1})$ or else a copy of the module $X_{\ell - 1}$ as
introduced in the proof of Lemma 13.  But as we argued before, this
contradicts finiteness of $\pdim N$ since $\ell - 1 \ge 0$. \qed 
\enddemo

\proclaim{Lemma 15} Let $\ell \ge 2$, and suppose that $N$ is an
indecomposable non-projective
$\la _{\ell}$-module of finite projective dimension such that
$a_{\ell}N \ne 0$ or $b_{\ell}N \ne 0$.  Then $N/JN$
contains a copy of $\bigl(S(b_\ell) \bigr)^m$.  In particular, $N$
requires at least $m$ generators.
\endproclaim

\demo{Proof} Invoking Lemma 14, we obtain $a_{\ell - 1}
\Omega^1(N)
\allowmathbreak \ne 0$.  Inspection of the radicals of the indecomposable
projective $\la _{\ell}$-modules now yields a top element of
type $a_{\ell}$ in $N$ which is annihilated by $\alpha _{\ell 0}$. 
Therefore our assertion follows from Lemma 12.
\qed \enddemo   

\proclaim{Lemma 16} All finitistic dimensions of $\la _0$ are equal to
$r$, that is,  
$$\findim_k \la_0 = \findim\la _0 = \Findim \la_0 =  r  \qquad \text{for
all\ } k.$$ 

\endproclaim

\demo{Proof} The inequality $\findim_1\la_0 \ge r$ is due to the fact
that $\pdim S(a_0) = r$.  To establish the inequality $\Findim\la _0 \le
r$, one checks that the {\bf s}-invariant of $\la _0$ (as in Theorem 1)
is $\pdim \la_0 \gamma _1 = \pdim S(c_2) = r - 2$, and then applies
Theorem 1.  \qed
\enddemo 

\subhead{Final claim of Part I}\endsubhead 
For the present choice of the
target function
$f$, the algebra $\la = \la _s$ satisfies all conditions listed in
Theorem 10.

\demo{Proof of the final claim} By construction, $J^3 = 0$.  To see that
$\la$ is even an $\left( \lceil s/2 \rceil + 1\right)$-stack of monomial
algebras, consider the following alternate stacking partition of the set
$E$ of vertices of $\la$: Namely, $E_0 = \Q _0^{(0)}$, and $E_\ell =
\left( \Q _0^{(2\ell - 1)} \cup \Q _0^{(2\ell)}\right)\setminus E_{\ell -
1}
$, whenever $\ell \ge 1$ and $2\ell \le s$.  If $s$ is odd, we define
$E_{\lceil s/2 \rceil}$ to be $\{a_s, b_s\} = \Q _0^{(s)}\setminus
E_{(s-1)/2}$.

We now establish the equalities $(\dagger )$ preceding the construction
of the $\la_\ell$.

Returning to the notation employed in the construction of the
algebras $\la_\ell$, we combine Lemma 16 with the last
statement of Lemma 11 to obtain, via an obvious induction on $\ell$, the
following family of inequalities:
$$r = \findim_k \la_0 \le \findim_k \la_\ell \le \findim \la_\ell \le
\findim \la_0 + \ell = r+\ell,$$
for all $k \in \NN$ and $\ell \le s$.

 Next we verify that $\findim_m\la
_\ell
\ge r +
\ell$, which, in view of the above inequalities, will show that
$\findim_k\la _\ell = r +
\ell$ for all $k \ge m$.  For each $\ell \ge 0$, we define a $\la
_\ell$-module
$N_\ell = P_\ell /V_\ell$ as follows: 
$$P_\ell = \la _\ell x_{\ell 0}\oplus
\operatornamewithlimits{\bigoplus}_{i = 1}^m \la _\ell x_{\ell i}$$ 
with
$x_{\ell 0} = a_\ell$ and $x_{\ell i} = b_\ell$ for $1 \le i \le m$,
i\.e\., $P_\ell = \la a_\ell \oplus (\la b_\ell)^m$, and $V_\ell
\subseteq J_\ell P_\ell$ is generated by $\gamma_0 x_{\ell 0}$ and
the differences $\alpha _{\ell i}x_{\ell 0} - \beta _\ell x_{\ell i}$ for
$1 \le i \le m$.  Note that $N_\ell$ has Loewy length 2 and a tree graph,
namely

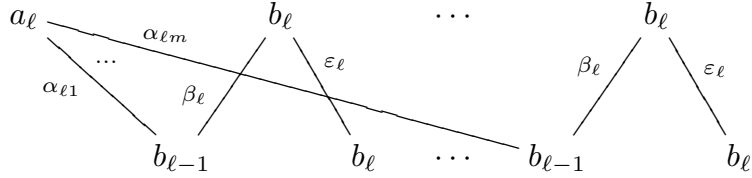
\ignore
$$\xymatrixcolsep{1pc}\xymatrixrowsep{2.5pc}
\xymatrix{ 
a_\ell \edge[drr]_(0.4){\alpha _{\ell 1}}^(0.4){\dots}
\edge[drrrrrr]^(0.25){\alpha_{\ell m}} &&& b_\ell
\edge[dl]_(0.65){\beta_\ell} \edge[dr]^(0.4){\varepsilon _\ell} &&
{\displaystyle{\cdots}} && b_\ell
\edge[dl]_(0.45){\beta _\ell} \edge[dr]^(0.45){\varepsilon _\ell} \\
 && b_{\ell-1} && b_\ell & {\displaystyle{\cdots}} & b_{\ell-1} &&
b_\ell
 }$$ 
\endignore

\noindent relative to
the obvious choice of top elements.  Moreover, we see that, for $\ell \ge
1$, we have $\Omega ^1(N_\ell) = N_{\ell - 1}$, while $\Omega ^1(N_0) =
S(c_1) \oplus (\la _0 b_{-1})^m$ has projective dimension $r - 1$ (see
\cite{\pre} for methodology).  Hence $\pdim N_\ell = r +
\ell$.

So only $\findim_{m-1}\la _{\ell}\le r$ for $1 \le \ell \le s$ remains to
be checked.  In light of Lemma 15, it suffices to verify this for $\ell =
1$.  Indeed, if $\ell \ge 2$ and $N$ is a module over $\la_\ell$, but not
over $\la_{\ell - 1}$, then  $N$ fails to be annihilated by at least
one of $a_\ell$, $b_\ell$. 

We write $\Delta = \la _1$ for ease of
notation and consider another stacking partition of the vertex set $E$
of $\Delta$.  Namely $E' = \{b_{-1}\}$ and
$E'' = E\setminus
E'$.  With $\be '$ and $\be ''$ as in Section 3, we define
$\Delta ' = \be '\Delta \be '$ and $\Delta '' =\be ''\Delta \be ''$ and
note that both $\Delta'$ and $\Delta''$ are monomial algebras.  It is
readily seen that $\findim \Delta' = 0$.  Moreover, we observe that
$\gamma _0$ is the only path
$p \in \eusm P''$ (the set of critical paths of the monomial algebra
$\Delta ''$) such that the maximum
$\text{\bf s}''$ of Theorem 1 is attained on $\Delta ''p$. In particular,
$\text{\bf s}'' = \pdim \Delta ''\gamma _0 = r - 1 \ge 1$.  Finally, we
note that all paths in $\eusm P''$ end in homogeneous vertices.

If $N$ is any $\Delta$-module with $\pdim_\Delta N = r + 1$, 
Corollary 9 therefore tells us that the $\Delta''$-projective dimension
of  
$\Omega ^2 _{\Delta''}(\be'' N)$ also equals $\text{\bf s}''$.
Consequently, Theorem 1 forces  $\Omega^2 _{\Delta''}(\be'' N)$ to
have a direct summand isomorphic to $\Delta'' \gamma_0$.  By the
second part of Theorem 1, we further obtain a top element $y \in
\be''N$ of type
$a_1$ with $\alpha _{10}y = 0$.  Clearly, $y$ is
also a top element of $N$, whence $N$ requires at least $m$ generators by
Lemma 12.  This proves $\findim_{m - 1}\Delta \le
r$ as required. 

That the final condition of Theorem 10 is met, is clear from the above
considerations.

\enddemo 

\subhead Part II of the proof of Theorem 10 \endsubhead  
Fix $m, \, n,
\, r, \, s, \, t \in \Bbb{N}$ with $2\le m < n$ and $r \ge 2,$ and assume
the target function $f$ of Theorem 10 to be 
$$\alignat2
                f(k) &= r &&\qquad \text{for} \quad 1 \le k \le m - 1, \\
                f(k) &= r + s &&\qquad \text{for} \quad m \le k \le n-1,
\\
\text{and} \quad f(k) &= r + s + t &&\qquad \text{for} \quad k \ge n.
\endalignat$$

For the construction of a finite dimensional algebra $\la$ having the
finitistic dimensions prescribed by $f$, let $\la _0, \la _1, \dots, \la
_{s-1}$ be as in Part I.  We slightly modify the final algebra $\la _s$
constructed before - its finitistic dimensions will remain unchanged - to
`switch gear' in our recursive pattern so as to smooth the road for
another jump.

\remark{Step $s$} Given $\la _0, \dots, \la _{s-1}$ as before and keeping
the notation pertaining to these algebras, we enlarge the quiver $\Q
^{(s-1)}$ of $\la _{s-1}$ to $\Q ^{(s)}$ as follows: We add the vertices
$a_s$, $b_s$, $b_{-1}'$, $b_0'$ and supplement arrows and
relations as indicated by the following graphs of the indecomposable
projective $\la _s$-modules:

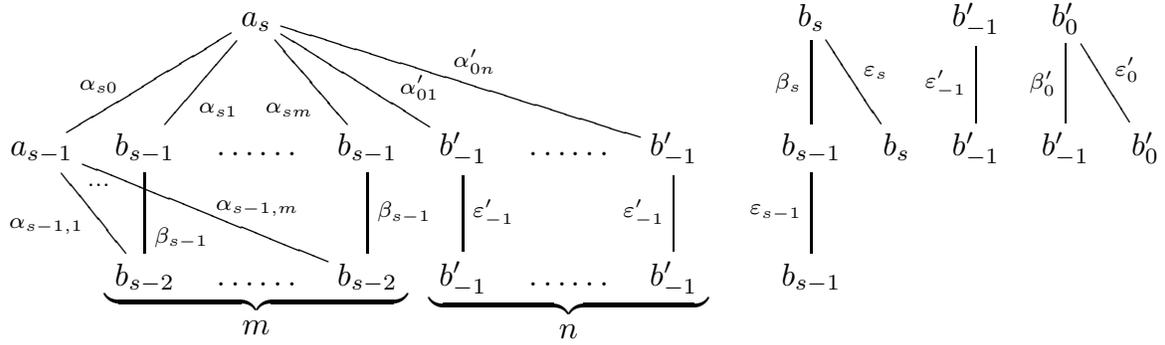
\ignore
$$\xymatrixcolsep{0.6pc}
\xymatrix{ & & a_s \edge[dll]_(0.65){\alpha _{s0}}
\edge[dl]^(0.6){\alpha _{s1}} \edge [dr]_(0.6){\alpha _{sm}} \edge
[drr]^(0.7){\alpha _{01}'} \edge [drrrr]^{\alpha _{0n}'} & & & & & & b_s
\edge[d]_{\beta_s} \edge[dr]^{\varepsilon _s} & & b_{-1}'
\edge[d]_{\varepsilon _{-1}'} & b_0'
\edge[d]_{\beta_0'} \edge[dr]^{\varepsilon _0'} \\
  a_{s-1} \edge[dr]_{\alpha _{s-1,1}}^(0.35){\dots}
\edge[drrr]^(0.6){\alpha _{s-1,m}} & b_{s-1}
\edge[d]^(0.7){\beta_{s-1}} &
\dots\dots & b_{s-1} \edge[d]^{\beta_{s-1}} & b_{-1}'
\edge[d]^{\varepsilon _{-1}'} &
\dots\dots & b_{-1}' \edge[d]_{\varepsilon _{-1}'} & & b_{s-1}
\edge[d]_{\varepsilon _{s-1}} & b_s & b_{-1}' & b_{-1}' & b_0' \\
 & b_{s-2} & \dots\dots \save[0,-1].[0,1]!C *\frm{_\}}="lbrc" \restore
\save"lbrc"+D \drop++!U{m} \restore
& b_{s-2} & b_{-1}' & \dots\dots \save[0,-1].[0,1]!C *\frm{_\}}="rbrc"
\restore \save"rbrc"+D \drop++!U{n}\restore
& b_{-1}' & & b_{s-1} }$$
\endignore

Here we preserve our convention concerning non-monomial relations of $\la
_s$ communicated by the graph of $\la _s a_s$; namely, that each subgraph
of type
$\tilde A_3$ of the graph of $\la_s a_s$ corresponds to a relation
$\alpha _{s-1, i}\alpha _{s 0} - \beta _{s-1}\alpha_{s i}$ in the ideal
$I$ of relations.

As before, $\findim_k \la_s \ge r$ for all $k$.  To check that
$\findim_k \la_s \le r$ for $k \le m-1$, let $N$ be a finitely generated
indecomposable
$\la_s$-module of projective dimension $r+1$.  If $s = 1$, the argument
given under the final claim of Step I shows $N$ to
require at least $m$ generators.  So suppose that $s \ge 2$.  If we can
show that $a_s N \ne 0$ or  $b_s N \ne 0$, then again the reasoning of
Step I provides what we need (Lemmas 14 and 15 carry over to the new
format of $\la_s$, with minor modifications of the arguments).  So it
suffices to consider the case where  $a_s N = b_s N = 0$.  Let $\be$ be
the sum of the primed primitive idempotents, that is, $\be = b'_{-1} +
b'_0$.  Our annihilation assumption ensures a direct-sum decomposition
of $N$ into
$\la_s$-modules $N_1$ and $N_2$ such that $\be N_1 = 0$ and $\be N_2 =
N_2$.  Since $N$ is indecomposable and $\pdim_{\la_s} N_2 = \pdim_{\be
(\la_s) \be} N_2 \le \findim \be (\la_s) \be = 0$, we infer that $N =
N_1$ is a $\la_{s-1}$-module.  In view of the known equality
$\findim_{m-1} \la_{s-1}=r$, we thus again glean a minimum of $m$
generators for $N$.       

To see that $\findim_k\la _s =
r + s$ for $k \ge m$, we can rely on the previous arguments.  However, the
presentation of the $m$-generated module
$N_s = P_s/V_s$ that has the same graph as the module of that name in Part
I needs to take the slightly altered structure of $\la _s$
into account as follows:  Again,
$$P_s = \operatornamewithlimits{\bigoplus}_{i=0}^{m}\la_s x_{si} \quad
\text{with } x_{s0} = a_s \text{ and } x_{si} = b_s \text{ for } 1 \le i
\le m,$$ but now $V_s$ is the submodule generated by $\alpha _{s0}
x_{s0}, \,\, \alpha _{si} x_{s0} - \beta _s x_{si}$, for $1 \le i \le m$,
\underbar{and} $\alpha _{0j}'x_{s0}$ for $1 \le j \le n$.  Clearly
$\Omega ^1(N_s) = N_{s-1} \oplus \bigl(\la_s b'_{-1}\bigr)^n$, whence
$\pdim N_s = r + s$. 

\endremark

\remark{Step $s+\ell$} for $\ell \ge 1$.  We introduce three new
vertices,
$a_{s+\ell}$, $b_{s+\ell}$, and
$b'_\ell$, giving rise to indecomposable projective left $\la
_{s+\ell}$-modules whose graphs revert to the mold of Part I.

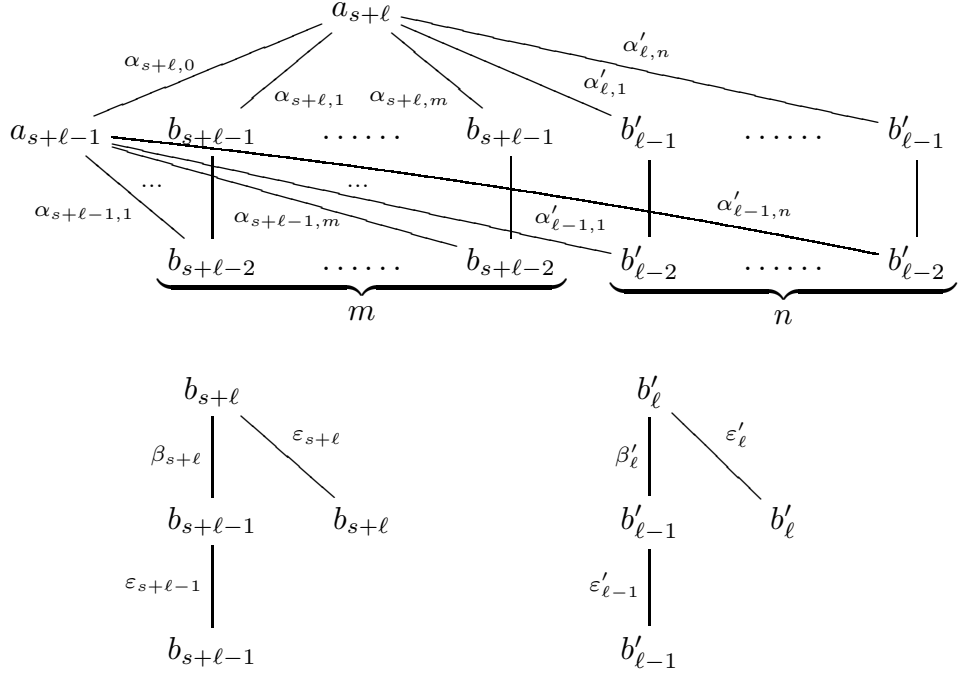
\ignore
$$\xymatrixcolsep{1.2pc}
\xymatrix{ & &  a_{s+\ell} \edge[dll]_(0.6){\alpha_{s+\ell,0}}
\edge[dl]^(0.6){\alpha_{s+\ell,1}}
\edge[dr]_(0.6){\alpha_{s+\ell,m}}
\edge [drr]^(0.8){\alpha'_{\ell,1}} \edge
[drrrr]^{\alpha'_{\ell,n}} \\
a_{s+\ell-1} \edge[dr]_{\alpha_{s+\ell-1,1}}^{\dots}
\edge[drrr]_(0.55){\alpha_{s+\ell-1,m}}
\edge[drrrr]^(0.85){\alpha'_{\ell-1,1}}^{\dots}
\edge@/^0.3pc/[drrrrrr]^(0.8){\alpha'_{\ell-1,n}} & b_{s+\ell-1}
\edge[d] &
\dots\dots & b_{s+\ell-1} \edge[d] & b'_{\ell-1} \edge[d]&
\dots\dots & b'_{\ell-1}
\edge[d] \\
 & b_{s+\ell-2} & \dots\dots \save[0,-1].[0,1]!C *\frm{_\}}="lbrc"
\restore  \save"lbrc"+D \drop++!U{m}\restore
& b_{s+\ell-2} & b'_{\ell-2} &
\dots\dots \save[0,-1].[0,1]!C *\frm{_\}}="rbrc"
\restore  \save"rbrc"+D \drop++!U{n}\restore
& b'_{\ell-2} \\
 & b_{s+\ell}
\edge[d]_{\beta_{s+\ell}} \edge[dr]^{\varepsilon_{s+\ell}} & & &
b'_{\ell} \edge[d]_{\beta'_{\ell}}
\edge[dr]^{\varepsilon'_{\ell}} \\
 & b_{s+\ell-1} \edge[d]_{\varepsilon_{s+\ell-1}} &
b_{s+\ell} & & b'_{\ell-1}
\edge[d]_{\varepsilon'_{\ell-1}} & b'_{\ell} \\
 & b_{s+\ell-1} & & & b'_{\ell-1} }$$
\endignore

\endremark

Lemmas 11-15 of Part I have analogues applying to the present situation. 
We list them in a format that permits us to carry over the previous
arguments almost verbatim.

\proclaim{Lemma 11$'$}  Let $\ell \ge 0$, and suppose that $N$ is a $\la
_{s+\ell +1}$-module of finite projective dimension.  Then $b_{s+\ell
+1}\Omega^1(N) = b'_{\ell
+1}\Omega^1(N) = 0$; in particular, $\Omega^1(N)$ is a $\la
_{s+\ell}$-module. \qed
\endproclaim

\proclaim{Lemma 12$'$}  Let $\ell \ge 0$, and suppose $N$ is a $\la _{s+\ell
+1}$-module of finite projective dimension containing a top element $y$
of type $a_{s+ \ell +1}$ such that $\alpha _{s+ \ell +1, 0} y = 0$.  Then
$N/JN$ contains $\bigl(S(b_{s+\ell +1})\bigr)^m$  and $\bigl(S(b'_{\ell
+1})\bigr)^n$. \qed  
\endproclaim

\proclaim{Lemma 13$'$} Let $\ell \ge 1$, and suppose $M$ is a $\la
_{s+\ell}$-module of Loewy length $2$ with $a_{s+\ell} M = 0$, but
$b_{s+ \ell} M \ne 0$ or $b'_\ell M \ne 0$.  Then $\pdim M =
\infty$.
\qed
\endproclaim 

To prove the case $\ell=0$ in the next lemma, note that, whenever $\pdim
N <\infty$ and $b'_1 N\ne 0$, the inequality $a_{s+1} N \ne 0$ is
automatic; hence, the hypothesis boils down to `$a_{s+1} N \ne 0$ or
$b_{s+1} N \ne 0$' in that case.

\proclaim{Lemma 14$'$} Let $\ell \ge 0$, and suppose that $N$ is an
indecomposable non-projective $\la _{s + \ell + 1}$-module of finite
projective dimension such that
$a_{s + \ell + 1}N \ne 0$ or
$b_{s + \ell + 1}N \ne 0$, or else $b'_{\ell+1}N \ne 0$.  Then
$a_{s+\ell}
\Omega ^1(N) \ne 0$. \qed
\endproclaim

\proclaim{Lemma 15$'$} {\rm(a)} Let $\ell \ge 1$, and suppose that $N$ is
an indecomposable non-projective 
$\la _{s+\ell}$-module of finite projective dimension with
$a_{s+ \ell}N \ne 0$ or $b_{s+\ell}N \ne 0$.  Then $N/JN$
contains a copy of $\bigl(S(b_{s+\ell}) \bigr)^m$.  In particular, $N$
requires at least $m$ generators.

{\rm(b)} Let $\ell \ge 2$, and suppose that $N$ is
an indecomposable  non-projective
$\la _{s+\ell}$-module of finite projective dimension with
$a_{s+ \ell}N \ne 0$ or $b_{s+\ell}N \ne 0$ or $b'_{\ell} N\ne 0$.  Then
$N/JN$ contains a copy of $\bigl(S(b'_\ell) \bigr)^n$.  In particular,
$N$ requires at least $n$ generators. \qed
\endproclaim

\subhead{Final claim of Part II} \endsubhead The algebra $\la = \la
_{s+t}$ has finitistic dimensions $\findim_k\la = f(k)$ for $k\in \NN$ and
satisfies the additional conditions of Theorem 10.

\demo{Proof of the final claim}  Here $d = s + t$.  To see that $\la$
is a
$\left( \lceil d/2 \rceil + 1 \right)$-stack of monomial algebras, we
proceed as in Part I, by setting $E_0 = \Q _0^{(0)}$ and adding on the
additional vertices of any pair $\la _{2\ell - 1}$, $\la _{2\ell}$ for
$2\ell \le s + t$ in the following steps to move from $E_{\ell - 1}$ to $E
_{\ell}$.  Since every left $\la_s$-module is also a
$\la_{s+\ell}$-module, we see that
$$r = \findim_k\la _s \le \findim_k\la _{s + \ell}$$ 
for $k \le m - 1$, and as before we obtain 
$$r + s = \findim_k\la _s \le \findim_k\la _{s + \ell} \le \findim \la _s
+ \ell = r + s + \ell$$ for $k \ge m$.  So we only need to verify the
following: \roster
\item $\findim_k\la _{s+\ell} \le r$ \quad for $k \le m - 1$ and all $\ell
\ge 1$
\item $\findim_k\la _{s+\ell} \le r + s$ \quad for $m \le k \le n - 1$ and
all
$\ell \ge 1$
\item $\findim_k\la _{s+\ell} \ge r + s + \ell $ \quad for $k \ge n$ and
all
$\ell \ge 1$.
\endroster

Concerning \therosteritem1: In view of Lemma 15$'$, it
suffices to prove $\findim_{m-1} \la_{s+1} \le r$.  We write $\Delta =
\la_{s+1}$ for convenience.  Let $N \in
\Delta\operatorname{-mod}$ be indecomposable of projective
dimension at least $r+1$.   If $a_{s+1} N \ne 0$ or $b_{s+1} N \ne 0$,
then $N$ requires at least $m$ generators by Lemma 15$'$(a). 
 So suppose that
$a_{s+1} N = b_{s+1} N = 0$.  But in that case, indecomposability of $N$
forces $b'_1$ to annihilate $N$ as well.  Indeed, since $N \not\cong
\Delta b'_1$, the inequality
$b'_1 N \ne 0$ would place a direct summand isomorphic to a  
submodule $X$ of $J b'_1$
 into the syzygy $\Omega^1(N)$; but this is not permissible because all
such modules $X$ have infinite projective dimension.  This makes
$N$ a $\la_s$-module, and guarantees that $N$ is not $(m-1)$-generated. 
The argument for \therosteritem1 is thus complete.

Concerning \therosteritem2:  Once more, Lemma 15$'$ restricts our
focus to $\Delta = \la_{s+1}$.  So suppose $N \in
\Delta\operatorname{-mod}$ is indecomposable  with $\pdim N \ge r+s+1$;
in particular, this implies $4 \le \pdim N < \infty$.  We aim at an
application of Lemma 12$'$ to show that
$N$ requires at least
$n$ generators.  In other words, we wish to show that the multiplicity
$\mu$ of $S(a_{s+1})$ in $N/JN$ exceeds the multiplicity $\nu$ of
$S(a_s)$ in
$JN/J^2N$.  Clearly $\nu \le \mu$. To obtain strict inequality, let $\be$
be the sum of all primitive idempotents of the form
$b_i$ and $b'_j$ in $\Delta$, and observe that  $\be \Delta \be$ is a
monomial algebra of finitistic dimension $\le 1$ (Theorem 1); in
particular, this ensures that neither $N$ nor $\Omega^1(N)$ is an $\be
\Delta
\be$-module.  Moreover, we observe that  
$\pdim X =
\infty$ for any indecomposable $\be \Delta \be$-module $X$ of Loewy
length at most $2$, except for $X = \Delta b_{-1}$, $\Delta b'_{-1}$,
$\Delta b_0$, and $\Delta b'_0$. 

First assume that $\mu = 0$, meaning
$a_{s+1} N = 0$.   In view of the properties of $\be
\Delta \be$, indecomposability of
$N$ then implies
$b_{s+1} N = b'_1 N = 0$, which makes $N$ a $\la_s$-module.  But this
contradicts $\findim \la_s = r+s$, and we conclude $\mu \ge 1$.  Once
more, we use indecomposability of $N$, combined with finiteness of $\pdim
N$, to see that $a_k(N/JN) = 0$ for all $k \le s$.  Consequently, the only
simple module of the form $S(a_j)$ potentially occurring in
$JN/J^2N$ is $S(a_s)$.  Since $a_k(J^2N) = 0$ for all $k$, equality $\nu
= \mu$ would therefore make $\Omega^1(N)$ an $\be \Delta \be$-module. 
This contradiction shows $\nu < \mu$, and an application of Lemma 12$'$
completes the proof of (2).

Concerning \therosteritem3: 
  For arbitrary choice of $\ell
\ge 1$, we consider the $n$-generated $\la _{s+\ell}$-mod\-ule 
$N_{s + \ell} = P_{s+\ell}/V_{s+\ell}$, where
$$P_{s+\ell} = \la_{s+\ell}\, x_{s+\ell,0} \oplus 
\biggl(\operatornamewithlimits{\bigoplus}_{i=1}^{m}\la_{s+\ell}\,
x_{s+\ell,i}\biggr) \oplus
\biggl(\operatornamewithlimits{\bigoplus}_{i=1}^{n}\la_{s+\ell}\,
x'_{\ell,i} \biggr);$$ 
here $x_{s+\ell,0} = a_{s+\ell}$, $x_{s+\ell,i} =
b_{s+\ell}$ for $1 \le i \le m$, and $x'_{\ell,i} = b'_{\ell}$ for $1 \le
i \le n$.  In other words,
$P_{s+\ell} = \la_{s+\ell}\, a_{s+ \ell} \oplus \bigl(\la_{s+\ell}\, b_{s+
\ell} \bigr)^m \oplus \bigl(\la_{s+\ell}\, b'_{\ell} \bigr)^n$.  Moreover
$V_{s+ \ell}$ is the submodule of $P_{s+ \ell}$ generated by
$\alpha_{s+\ell,0} x_{s+\ell,0}$, the differences $\alpha_{s+ \ell,i}
x_{s+ \ell,0} -  \beta_{s + \ell} x_{s+\ell,i}$ for $1 \le i \le m$, and
the differences $\alpha'_{\ell,i} x_{s+ \ell,0} -  \beta'_{\ell}
x'_{\ell,i}$ for $1 \le i \le n$.  Then $N_{s+ \ell}$ has the following
graph relative to the listed top elements:

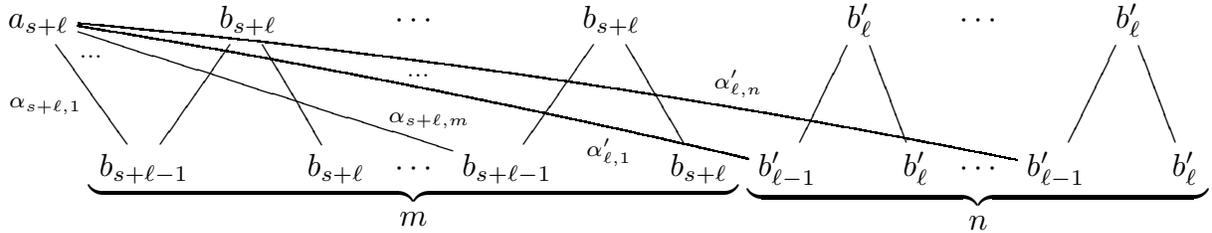
\ignore
$$\xymatrixcolsep{0.3pc}\xymatrixrowsep{2.5pc}
\xymatrix{ a_{s+\ell} \edge[dr]_{\alpha_{s+\ell,1}}^(0.3){\dots}
\edge[drrrrr]^(0.8){\alpha_{s+\ell,m}}
\edge@/^0.2pc/[drrrrrrrr]_(0.77){\alpha_{\ell,1}'}^{\dots}
\edge@/^0.4pc/[drrrrrrrrrrrr]^(0.68){\alpha_{\ell,n}'} & & b_{s+\ell}
\edge[dl]
\edge[dr] & & {\displaystyle{\cdots}} & & b_{s+ \ell}
\edge[dl] \edge[dr] & & & b_{\ell}' \edge[dl] \edge[dr] & &
{\displaystyle{\cdots}} & & b_{\ell}' \edge[dl] \edge[dr] \\
 & b_{s+\ell-1} & & b_{s+\ell} & {\displaystyle{\cdots}} 
\save[0,-3].[0,3]!C *\frm{_\}}="lbrc"
\restore  \save"lbrc"+D \drop++!U{m}\restore &
b_{s+\ell-1} & & b_{s+\ell} & b_{\ell-1}' & & b'_{\ell}
& {\displaystyle{\cdots}} \save[0,-3].[0,3]!C *\frm{_\}}="rbrc"
\restore  \save"rbrc"+D \drop++!U{n}\restore
& b_{\ell-1}' & & b'_{\ell} }$$ 
\endignore

\noindent As in Part I, we have $\Omega^1(N_{s + \ell}) = N_{s + \ell
-1}$ whenever $\ell
\ge 2$; moreover, $\Omega^2(N_{s+1}) = \Omega^1(N_s)$.  In light of
$\pdim N_s = r + s$, we conclude that
$\pdim N_{s + \ell} = r + s + \ell$ as desired. \qed\enddemo

\enddocument